\documentclass[12pt]{amsart}
\usepackage[margin=1in]{geometry}
\usepackage[utf8]{inputenc}
\usepackage{amsmath}
\usepackage{amssymb}

\usepackage[pdftitle={DisquesPFH4},
  pdfauthor={Samuel Bronstein},
  pdffitwindow=true,
  breaklinks=true,
  colorlinks=true,
  urlcolor=blue,
  linkcolor=red,
  citecolor=red,
  anchorcolor=red]{hyperref}
\usepackage{titleps}
\usepackage{xcolor}
\definecolor{links}{RGB}{50,0,200}
\definecolor{hyperrefcolor}{rgb}{50,0,250}
\hypersetup{citecolor=blue,linkcolor=links}

\usepackage{subfiles}%For other tex files
\usepackage{pgfplots}
\pgfplotsset{compat=1.15}
\usepackage{mathrsfs}
\usetikzlibrary{arrows}
\tikzset{naming/.style={align=center,font=\footnotesize}}
\tikzset{area/.style = {draw, shape = regular polygon, regular polygon sides = 6, thick, minimum width = 5cm}}

\numberwithin{equation}{section}
\numberwithin{table}{section}
\numberwithin{figure}{section}

\newcommand{\D}{\mathbb D}
\newcommand{\R}{\mathbb R}
\renewcommand{\H}{\mathbb H}
\renewcommand{\S}{\mathbb S}

\newcommand{\eps}{\varepsilon}
\newcommand{\Vol}{\text{Vol}}
\newcommand{\avg}[1]{\overline{#1}}
\newcommand{\Supp}{\text{Supp}}
\newcommand{\one}{\mathbf{1}}

\theoremstyle{plain}
\newtheorem{theorem}{Theorem}[section]
\newtheorem*{theorem*}{Theorem}
\newtheorem{lemma}[theorem]{Lemma}
\newtheorem{proposition}[theorem]{Proposition}

\newtheorem{corollary}[theorem]{Corollary}
\newtheorem*{corollary*}{Corollary}

\theoremstyle{definition}
\newtheorem{definition}[theorem]{Definition}

\theoremstyle{remark}
\newtheorem{remark}{Remark}[section]

\title{On the Moser-Trudinger inequality for surfaces}
\author{samuel bronstein}
\date{2023}

\begin{document}
\maketitle
\begin{abstract}
	This paper is devoted to the Moser-Trudinger inequality on smooth riemannian surfaces.
	We establish that the constants involved can be chosen to depend on only $3$ parameters,
	which are the systole, isoperimetric constant and curvature of the surface.
	We have two analogous statements, considering respectively infinite-volume surfaces
	and closed surfaces.
	As an application, we show that there are sequences of coverings of a hyperbolic closed
	surface which admit a uniform Moser-Trudinger inequality.
\end{abstract}
\tableofcontents
\section{Introduction}
In this paper we consider $(\Sigma,\sigma)$ a Riemannian surface,
and we look at the famous Sobolev embeddings $W^{1,2}(\Sigma)\hookrightarrow L^p(\Sigma)$
for any $p\geq 2$.

When $\Sigma$ is of infinite volume, we prove the following version of the celebrated
Moser-Trudinger inequality:
\begin{theorem}
	Let $(\Sigma,\sigma)$ be an infinite volume complete surface.
	Assume that the curvature of $\Sigma$ is upper bounded by $K$.
	Assume that the systole $\delta$ and the Cheeger isoperimetric constant $h$ are
	both nonzero.
	Then there is a constant $C$ depending only on $(\delta,h,K)$ such that, for any
	$u\in W^{1,2}(\Sigma)$:
	\begin{equation}
		\int_\Sigma|\nabla u|^2 d\sigma\leq 1\,\Rightarrow\,
	\int_\Sigma\exp\big(4\pi u^2\big)-1 d\sigma\leq C(\delta,h,K)
	\end{equation}
\end{theorem}
We refer to section $2$ for precise definitions of systole and isoperimetric constant for
noncompact surfaces.

For compact surfaces, this statement cannot be true because of the existence of nonzero harmonic
maps.
The standard way to deal with this problem is to add the assumption that our maps have zero average,
restricting effectively to the orthogonal of harmonic maps in $W^{1,2}$.

In this setting, we prove the following inequality:
\begin{theorem}
	Let $a>0$ and $K\in\R$. There is a constant $C(a,\delta,h,K)$ such that,
	for any closed surface $(\Sigma,\sigma)$ with systole $\delta$,
	Cheeger constant $h$, curvature upper bounded by $K$ and volume bigger than $a$,
	for any $u\in W^{1,2}(\Sigma)$ satisfying:
	\begin{equation}
		\int_\Sigma u d\sigma =0,\quad\int_\Sigma |\nabla u|^2d\sigma\leq 1
	\end{equation}
	Then:
	\begin{equation}
		\int_\Sigma \exp(4\pi u^2)-1 \,d\sigma\leq C(a,\delta,h,K)
	\end{equation}
\end{theorem}

When the considered surface is the sphere $\S^2$ with its round metric, our theorem is
exactly the Moser-Trudinger inequality, first proved by Trudinger \cite{Tru67}
and then generalized by Moser to the higher-dimensional spheres \cite{Mos71}.

When considering the hyperbolic plane $\H^2$ which satisfies the condition of the abovementioned theorem,
we get back a result of Mancini and Sandeep \cite{MS10}.

The constant $4\pi$ is the biggest one for which the statement is true,
statement proven by Li \cite{Li01} for any compact Riemannian surface.

We insist that the main novelty is the dependency of the constant $C$ involved in the parameters, and
particularly that it is roughly independent of the volume of the surface.
There are other ways to prove the Trudinger-Moser inequality for compact surfaces, by a concentration argument, as Li did \cite{Li01}.
We would expect other more recent methods to work to, for instance the fractional Sobolev
inequalities have recently been used to prove an analogous inequality in the case of the $n$-dimensional sphere , see Xiong \cite{Xio18}.

The same kind of dependency is also true for the constants involved in the Sobolev embeddings,
see Ilias's work \cite{Ili83}. See also Faget \cite{Fag06}, for a generalization of the exponential
decay of Sobolev's constants in higher dimensional manifolds.

This reinforces the feeling that Moser-Trudinger is a nice way of understanding the asymptotic
behavior of the embeddings $W^{1,2}(\Sigma)\hookrightarrow L^p(\Sigma)$.

Generally, we have no idea of the sharp constant $C$ involved in the theorem.
This question was asked in a paper of Dolbeault-Esteban--Jankowiak \cite{DEJ14}.
Our paper gives a partial answer as it claims that the constant won't explode if we control
the right parameters, but we think the exhibition of the sharp constant should involve more
complicated techniques.

Our proof mainly relies on a symmetrisation argument, which naturally involves the isoperimetric
inequality. That's why we preferred to write the dependency in terms of the isoperimetric constant
rather than the spectral gap, even if they are closely linked by the Cheeger and Buser inequalities \cite{Che70, Bus82}. 

Trudinger-Moser inequality has had tremendous consequences regarding the prescribed curvature
problem on $\S^2$. Famously, Moser used it to solve the problem of prescribed curvature
on the projective plane \cite{Mos73}.
Mancini--Sandeep \cite{MS10} use their version of Moser-Trudinger for the prescribed curvature
problemn on subsets of $\R^2$.
It has also been used by Jost--Wang \cite{JW01} to solve some Toda systems on
closed surfaces. Malchiodi--Ruiz \cite{MR11} and \cite{BJM+15} deduced existence theorems for
mean-field equations from the Trudinger-Moser inequality.

Our focus here is rather on hyperbolic surfaces of high genus.
A corollary of our theorem is the possibility of taking high degree covers of a fixed hyperbolic
surface keeping a uniform sobolev embedding
\begin{corollary}
	Let $\Sigma$ be a closed hyperbolic surface.
	There are sequence $(\Sigma_n)$ of Riemannian covers of $\Sigma$ of degree $n$
	and $C(\Sigma)>0$ such that, for any $n\geq 1$,
	$u\in W^{1,2}(\Sigma_n)$ satisfying
	\begin{equation}
		\int_{\Sigma_n}u=0,\quad\int_{\Sigma_n}|\nabla u|^2\leq 1
	\end{equation}
	then
	\begin{equation}
		\int_{\Sigma_n} e^{4\pi u^2}-1\leq C(\Sigma)
	\end{equation}
\end{corollary}
Indeed, thanks to Magee--Naud--Puder \cite{MNP22} we know that for any $\eps>0$,
with probability going to one when the index of the covering increases,
a random cover of $X$ has relative spectral gap bigger than $\frac{3}{16}-\eps$.
This means that there are sequences of coverings of a closed surface with spectral gap uniformly
lower bounded, hence for such a sequence $(X_n)$ the Cheeger constant is lower bounded.
The systole and volume are bigger than the systole and volume of $X$, and the curvature is $-1$,
hence the constant in Moser-Trudinger is uniform.

The author is very thankful to Nicolas Tholozan and Paul Laurain for their time and discussion on
these topics.

\section{Notations and isoperimetric comparison}
Throughout this section $(\Sigma,\sigma)$ will be a Riemannian surface.
Recall the definitions of systole and  Cheeger constant:
\begin{definition}[Systole]\label{sys}
The systole of $(\Sigma,\sigma)$ is the length of the shortest noncontractible closed curve
in $\Sigma$.
\end{definition}
\begin{definition}[Cheeger constant]\label{chect}
Denote $l$ the length measure associated to $\sigma$.
The Cheeger constant of $(\Sigma,\sigma)$ is:
\begin{equation}
	h:=\underset{A\subset\Sigma,\,2\Vol(A)<\Vol(M)}{\inf}\frac{l(\partial A)}{\Vol(A)}
\end{equation}
If $\Vol(\Sigma)=\infty$, we require $\overline A$ to be compact.
\end{definition}
The Cheeger constant is heavily linked with the spectral gap of the surface.
\begin{definition}
For a closed surface, the spectral gap $\Lambda$ of $(\Sigma,\sigma)$ is the smallest nonzero
eigenvalue of $-\Delta$.
\end{definition}
It is a standard fact that the spectrum of the Laplacian on a compact surface is discrete,
hence the spectral gap of a closed surface is always strictly positive.
Remark that there is a way to generalize the spectral gap for asymptotically negatively curved
surfaces, using the resolvent of the laplacian. See for instance Borthwick \cite{Bor07} for the
spectral gap of infinite area hyperbolic surfaces.
On the hyperbolic plane, the spectral gap will be $\frac{1}{4}$, which is a reformulation
of Hardy's inequality.
The spectral gap and Cheeger's isoperimetric constant are closely linkedn thanks to two relations.
Recall the first one, proved by Cheeger in \cite{Che70}:
\begin{theorem}[Cheeger]
Let $\Sigma$ be a complete surface. Denote $h$ and $\Lambda$ its Cheeger constant and spectral gap.
Then
\begin{equation}
	4\Lambda\geq h^2\,.
\end{equation}
\end{theorem}
Remark that Cheeger states it for closed surfaces, but it has since been extended for complete
noncompact surfaces, see \cite{BPP93}.
Buser \cite{Bus82} , theorem 7.1., showed that for surfaces with curvature bounded below,
a converse holds:
\begin{theorem}[Buser]
There is a universal constant $C$, such that:
Let $(\Sigma,\sigma)$ be a complete Riemannian surface with curvature bounded below by
$-\delta^2$, $\delta\geq 0$. Denote $h$ and $\Lambda$ its Cheeger constant and spectral gap. Then
\begin{equation}
	\Lambda\leq C(\delta h+h^2)
\end{equation}
\end{theorem}
We will need a stronger insight on the isoperimetric behavior of a surface, hence we introduce
the \emph{isoperimetric profile}, as considered by Grimaldi--Nardulli--Pansu \cite{GNP09}:
\begin{definition}[Isoperimetric profile]\label{isopro}
Let $(\Sigma,\sigma)$ be a Riemannian surface.
The isoperimetric profile is the following map:
\begin{equation}
	\forall t>0,\,\varphi(t)=\underset{A\subset\Sigma,\Vol(A)=t}{\inf} l(\partial A)
\end{equation}
\end{definition}
Remark that this notation differs from the isoperimetric profile for revolution manifolds
considered by Gallot in \cite{Gal88}, but it is easy to go from one definition to another.
\begin{lemma}\label{smallprofile}
Let $\Sigma$ be a complete surface with nonzero systole $\delta$ and
curvature bounded above by $K$. Denote $\varphi$ its isoperimetric profile.
Denote $V_0$ the following:
\begin{equation}
	V_0 =\frac{\delta^2}{2\pi+\sqrt{4\pi^2-K\delta^2}}
\end{equation}
	Consider $t_0=\min\{V_0,\frac{\Vol(\Sigma)}{2}\}$ if $K\leq 0$
	and $t_0=\min\{V_0,\frac{\Vol(\Sigma)}{2},\frac{4\pi}{1+K}\}$ if $K>0$
	Then for $t\leq t_0$, we have:
\begin{equation}
	\varphi(t)\geq \sqrt{t(4\pi-Kt)}
\end{equation}
\end{lemma}
\begin{proof}
Assume first $V_0\leq\frac{\Vol(\Sigma)}{2}$
Let $A\subset\Sigma$ be a domain of area $t\leq V_0$.
$V_0$ was chosen so that:
\begin{equation}
	\sqrt{V_0(4\pi-K V_0)}\leq \delta
\end{equation}
Hence the wanted inequality is obviously verified if $l(\partial A)\geq \delta$.

Otherwise, if $l(\partial A)\leq\delta$, because of completeness one can assume that $A$ is a
disjoint union of disks. But for a disk, as the curvature is lower than $K$,
the length of its boundary is bigger than the length of a disk of same area
in constant curvature $K$ model space.
Hence, we get the inequality desired.

\end{proof}
We will compare the isoperimetric profile of a surface with the behavior of a radial metric on
the disk.

\begin{lemma}\label{globalprofile}
	Consider $\eps>0$, $h>0$ and $K\in\mathbb R$.
	Let $g$ be defined as the following map:
\begin{equation}
	\forall t >0,\,g(t):=\left\{\begin{array}{ccc}
		\sqrt{t(4\pi-K t)}&\text{ if }& t\leq\eps\\
		g(\eps)&\text{ if}& \eps< t\leq \frac{g(\eps)}{h}\\
		h\cdot t&\text{ if }& t> \frac{g(\eps)}{h}
	\end{array}\right.
\end{equation}
Then there exists a $C^1$ map $f:\R_+\rightarrow \R_+$ such that the disk endowed with radial
metric $dr^2+f'(r)^2d\theta^2$ satisfies:
\begin{enumerate}
	\item
		The perimeter of a centered disk of radius $A$ is $g(A)$
	\item
		there is a constant $C(g)>0$ such that the metric $g_\D=dr^2+f'(r)^2d\theta^2$ is uniformly conformal to a disk with constant negative curvature.
\end{enumerate}
\end{lemma}
\begin{figure}[htb]
\centering
\input{radial.tex}
\caption{graph of g for $(K,h,\eps)=(-1,1.6,0.2)$}
\end{figure}
\begin{proof}
The assumption (1) is verified if and only if $f$ satisfies the functional equation:
\begin{equation}
	2\pi f'(r)=g(2\pi f(r))\,\forall r>0
\end{equation}
This equation can be explicitly solved on $(0,\eps)$
with:
\begin{equation}
	f(r)=\left\{\begin{array}{ccc}
		\frac{-2}{K}\sinh(\frac{\sqrt{-K}r}{2})&\text{ if }&K<0\\
		\frac{r^2}{2}&\text{ if }&K=0\\
		\frac{2}{K}\sin(\frac{\sqrt{K}r}{2})&\text{ if }&K>0
	\end{array}\right.
\end{equation}
On $(\eps,\infty)$ there is a nondecreasing solution on $\R$, which has at least linear growth,
as $g$ is lower bounded. The solution is defined on all of $\R$ because $g$ is sublinear.
Actually, we have the following explicit description of $f$:
\begin{flalign}
	\forall \eps< t\leq \frac{g(\eps)}{h},\,&f(t)=f(\eps)+(t-\eps)\frac{g(\eps)}{2\pi}\\
	\forall t>\frac{g(\eps)}{h},\,&f(t)=f(\frac{g(\eps)}{h})\exp\big(h(t-\frac{g(\eps)}{h})\big)
\end{flalign}
As for $t>\frac{g(\eps)}{h}$, we have:
\begin{equation}
\frac{f'(t)}{\sinh(ht)}=2hf(\frac{g(\eps)}{h})e^{-g(\eps)}\frac{1}{1-e^{-2ht}}
\end{equation}
And $\frac{f'(t)}{\sinh(ht)}$ has a nonzero limit when $t\rightarrow 0$.
We deduce that there is a constants $C(K,\eps,h)$ so that
\begin{equation}
\frac{1}{C(K,\eps,h)}\leq\frac{f'(t)}{\sinh(ht)}\leq C(K,\eps,h)
\end{equation}
so $f'$ is uniformly comparable to $\sinh(hr)$, which means that our disk is uniformly
quasiconformal to the complete disk with constant sectional curvature $-h$.
As this disk is quasiconformal to the hyperbolic plane, we get the statement (2).
We point out that the quasiconformal factor depends only on  the map $g$, so in 
the parameters $(K,\eps, h)$.
\end{proof}

\section{The Moser-Trudinger inequality for infinite-volume surfaces}
It is known that a Moser-Trudinger type inequality can't be true for all noncompact surfaces,
as it implies the existence of a Poincaré inequality.
We will show here that for a class of surfaces which at infinity, look like they are negatively
curved, a Moser-Trudinger inequality is true. We will also show that the relevant geometric data
are the systole and Cheeger constant. 
Throughout this section $(\Sigma,\sigma)$ will be a Riemannian surface, noncompact.
We will show that the following geometric assumptions are enough to get a Moser-Trudinger inequality:
\begin{theorem}\label{TMnoncpct}
Let $(\Sigma,\sigma)$ be an infinite volume complete surface.
Assume that the curvature of $\Sigma$ is upper bounded by $K$.
Assume that the systole $\delta$  and the Cheeger isoperimetric constant $h$ are both nonzero.
Then there is a constant $C$ depending only on $\delta,h,K$ such that,
for any $u\in W^{1,2}(\Sigma)$:
\begin{equation}
	\int_\Sigma|\nabla u|^2 d\sigma\leq 1\,\Rightarrow\,
	\int_\Sigma\exp\big(4\pi u^2\big)-1 d\sigma\leq C(\delta,h,K)
\end{equation}
\end{theorem}
Note that the Moser-Trudinger with exponent $4\pi$ was already shown by Mancini--Sandeep \cite{MS10}
in the case of disks dominated by the hyperbolic metric.
Our proof relies on a radial rearrangement adapted to the surface $\Sigma$,
the same used by Ilias \cite{Ili83} to estimate the constants involved in the Sobolev embeddings.
The following formulas are well-known measure-theoretic tools:
\begin{lemma}\label{formulas}
Let $u$ be a positive $C^2$ Morse function. Assume it is $W^{1,2}(\Sigma)$.
Denote $\nu$ the length measure associated to $\sigma$.
Consider the functions $l(t)$ and $A(t)$:
\begin{flalign}
	A(t)&:=\sigma\{u\geq t\}\\
	l(t)&:=\nu\{u=t\}
\end{flalign}
Then, for any measurable function $f:\R\rightarrow\R$ :
\begin{flalign}
	\int_\Sigma f(u) &=\int_\R f(t)A(t)dt\\
	\int_\Sigma|\nabla u|^2&\geq\int_\R \frac{l^2(t)}{A'(t)}dt
\end{flalign}
\end{lemma}
\begin{proof}
The first term is a standard formula in differential geometry, it is proven via a change of variables.
For the second inequality, we use the following formula for the derivative of $A$:
\begin{equation}
	A'(t)=\int_{\{u=t\}}\frac{d\nu}{|\nabla u|}
\end{equation}
Now the co-area formula yields (see \cite{Kes06}, theorem 2.2.2)
\begin{equation}
	\int_\Sigma |\nabla u|^2=\int_\R\big(\int_{\{u=s\}}|\nabla u|d\nu\big)ds
\end{equation}
And finally, applying Cauchy-Schwarz:
\begin{equation}
	\int_{\{u=s\}}|\nabla u|d\nu\int_{\{u=s\}}\frac{d\nu}{|\nabla u|}\geq l^2(s)
\end{equation}
which finishes the proof.
\end{proof}
\begin{theorem}\label{radialrearrangementnoncompact}
	Let $(\Sigma,\sigma)$ be an infinite volume complete surface,
	with nonzero systole, Cheeger constant and upper bounded curvature $\delta, h,K$.
	There is an nondecreasing function $f$ depending only on $(\delta,h,K)$,
	a disk uniformly conformal to the complete hyperbolic plane,
	with radial metric $dr^2+f'(r)^2d\theta^2$.
	Moreover, the radial rearrangement mapping on that disk:
	\begin{equation}
		\left\{\begin{array}{cc}
			W^{1,2}(\Sigma)&\rightarrow W^{1,2}(\D)\\
		u&\mapsto u^\ast\end{array}\right.
	\end{equation}
	defined so that $u^\ast$ is nonincreasing, radial, and:
	\begin{equation}
		A(t)=|\{u^\ast\geq t\}|\quad\forall t\in\R
	\end{equation}
	satisfies:
	\begin{flalign}
	 \int_\Sigma f(u)d\sigma&=\int_\D f(u^\ast)\\
	\int_\Sigma |\nabla u|^2&\geq \int_\D|\nabla u^\ast|^2
	\end{flalign}
\end{theorem}
\begin{proof}
	Denote $\varphi$ the isoperimetric profile of $\Sigma,\sigma$
	By lemma \ref{smallprofile}, there is $V_0(K,\delta)>0$ such that, on $(0,V_0)$
	the isoperimetric profile satisfies:
	\begin{equation}
		\varphi(t)\geq\sqrt{t(4\pi-Kt)}
	\end{equation}
	Consider the map $g$ defined in lemma \ref{globalprofile} with $\eps=V_0$.
	Then we claim:
	\begin{equation}
		\forall t \in\R,\,\varphi(t)\geq g(t)
	\end{equation}
	Indeed, for small values it is the lemmma \ref{smallprofile}, for big values it is the definition
	of the Cheeger constant and in between it comes from the fact that $g(V_0)\leq \delta$.

	Hence define $f$ by the lemma \ref{globalprofile}.
	For $u$ $C^2$, Morse, nonnegative and $W^{1,2}$ on $\Sigma$ define
	$u^\ast$ its radial rearrangement on $\D,dr^2+f'(r)^2d\theta^2$.
	By construction, it is clear that
	\begin{equation}
		\int_\Sigma f(u)=\int_\D f(u^\ast)
	\end{equation}
	for any measurable $f$.
	Also, denote $l$ and $l^\ast$
	the length of level set functions:
	\begin{flalign}
		l(t)&=\nu\{u=t\}|\\
		l^\ast(t)&=|\{u^\ast=t\}
	\end{flalign}
	Because $g\leq\varphi$ everywhere, it is clear that $l\geq l^\ast$.
	Thanks to lemma \ref{formulas}, we deduce
	\begin{equation}
	\int_\Sigma |\nabla u|^2\geq\int_\D|\nabla u^\ast|^2
	\end{equation}
\end{proof}
\begin{proof}[Proof of theorem \ref{TMnoncpct}]
	Consider $\Sigma$ such a surface.
	Consider $u\in W^{1,2}(\Sigma)$. First we will prove the inequality for Morse surfaces,
	then by density it will be true on $W^{1,2}(\Sigma)$.

	As the functional $e^{4\pi x^2}$ is even  and the absolute-value is $1$-lipschitz,
	up to replacing $u$ by $|u|$, we can assume that $u$ is nonnegative and
	\begin{equation}
		\int_\Sigma |\nabla u|^2\leq 1\,.
	\end{equation}
	Consider $(\D,dr^2+f'(r)^2d\theta^2)$ the disk defined in
	theorem~\ref{radialrearrangementnoncompact} and $u^\ast$ the radial rearrangement.
	By construction,
	\begin{flalign}
		\int_\D \exp\big(4\pi (u^{\ast})^2\big)-1&=\int_\Sigma \exp\big(4\pi u^2\big)-1\\
		\int_\D |\nabla u^\ast|^2&\leq 1
	\end{flalign}
	As $\int_\D |\nabla u^\ast|^2$ is a conformal invariant, consider $C$ the conformal
	distance between $dr^2+f'(r)^2d\theta^2$ and the metric of the Poincaré disk.
	Naturally, with $C$ being the conformal distortion, we control:
	\begin{equation}
		\int_0^\infty (\exp\big(4\pi(u^\ast)^2\big)-1)f'(r)dr\leq
		C\int_0^\infty (\exp\big(4\pi(u^\ast)^2\big)-1)\sinh(r)dr
	\end{equation}
	Now, cf \cite{MS10}, this is upperbounded by some bound depending only on $h$.
	Finally, as $C$ depends only on $(\delta, h,K)$ we get a bound depending on those
	parameters.
\end{proof}
We finish our discussion about infinite volume surfaces with $3$ examples of noncompact surfaces
which do not satisfy the Moser-Trudinger
inequality, showing that none of the $3$ assumptions could have been dropped.
\begin{theorem}
	Let $\Sigma$ be a (noncompact) surface with Cheeger constant $0$,
	and lower bounded curvature.
	Then the Moser-Trudinger inequality is false.
\end{theorem}
\begin{proof}
	Thanks to the work of Buser \cite{Bus82}, we know that such a surface has no spectral gap,
	hence it cannot satisfy a Moser-Trudinger inequality.
\end{proof}
\begin{remark}
	There are examples of surfaces from Buser \cite{Bus79} of surfaces with
	Cheeger constant arbitrarily small with regard to the spectral gap, but these examples
	involve small systoles and very negative curvature. It is not known to the author if with
	those examples one could show a noncompact surface satisfying the Moser-Trudinger inequality,
	but with Cheeger constant zero.
\end{remark}
\begin{theorem}
	Let $\Sigma,\sigma$ be a (noncompact) surface with curvature $\kappa$ satisfying:
	\begin{equation}
		\sup \kappa=\infty
	\end{equation}
	Then the Moser-Trudinger inequality is false.
\end{theorem}
\begin{proof}
	The main tool for this is to show that the constant of Moser-Trudinger explodes in high
	positive curvature.
	Fix $M>0$
	Let $\Omega$ be a small disk on which $\kappa\geq M$
	Consider $C_\S$ the Moser-Trudinger constant for the sphere.
	Thanks to works of Chang and Yang, there is a bubbling sequence $u_n$ of functions on $\S^2$,
	weakly converging towards a Dirac distribution,
	so that:
	\begin{flalign}
		\int_{\S^2}|\nabla u_n|^2&\leq 1\\
		\int_{\S^2}\exp\big(4\pi u_n^2)-1&\rightarrow C_\S
	\end{flalign}
	Because of the weak convergence, one can assume that $\Supp u_n$ is in small neighborhood of 
	a point.
	Hence for any $\eps>0$ we easily get a function $v$ supported in $\Omega$ satisfying:
	\begin{flalign}
		\int_\Sigma |\nabla v|^2&\leq 1\\
		\int_\Sigma \exp(4\pi v^2)-1\geq MC_\S-\eps
	\end{flalign}
	So the Moser-Trudinger inequality is false.
\end{proof}
When the systole goes to zero, there are plenty known examples where the inequality is false.
\begin{theorem}
	Let $\Sigma$ be the hyperbolic cusp:
	$\R\times\S^1$ equipped with the metric $dt^2+e^{2t}d\theta^2$.
	Then it does not satisfy the Moser-Trudinger inequality.
\end{theorem}
\begin{proof}
	Fix $b\in\R$, and consider the function
	\begin{equation}
		f(t,\theta)=(e^t-e^b)\one_{t\leq b}
	\end{equation}
	One can make the following computations:
	\begin{flalign}
		|\nabla f|_2^2&=\frac{2\pi}{3}e^{3b}\\
		|f|_2^2&=\frac{2\pi}{3}e^{3b}\\
		|f|_4^4&=\frac{2\pi}{5}e^{5b}
	\end{flalign}
	hence we see that, when $b\rightarrow-\infty$
	\begin{equation}
		\frac{|f|_4^4}{|\nabla f|_2^4}\rightarrow +\infty
	\end{equation}
	which contradicts the Moser-Trudinger inequality.
\end{proof}

\section{A Moser-Trudinger inequality for closed surfaces}
In this section, we need to adapt our proof with more subtlety in order to get the same theorem
for closed surface.
We point out that Li \cite{Li01} already proved that all compact surfaces
satisfy the Moser-Trudinger inequality with exponent $4\pi$, however
his proof does not show the dependency of the upper bound
in the systole, curvature and Cheeger constant.

\begin{theorem}\label{TMcompact}
Let $\Sigma,\sigma$ be a closed surface with nonzero systole $\delta$, Cheeger constant $h$
and upper bounded curvature $K$.
	Then there is a constant $C(\delta,h,K,\Vol(\Sigma))$ such that,
for any $u\in W^{1,2}(\Sigma)$ with zero average and satisfying:
\begin{equation}
\int_\Sigma |\nabla u|^2d\sigma\leq 1
\end{equation}
Then
\begin{equation}
\int_\Sigma\exp\big(4\pi u^2\big)-1\,d\sigma\leq C(\delta,h,K)\,.
\end{equation}
	The constant $C$ can be chosen to be nonincreasing in $\Vol(\Sigma)$.
\end{theorem}
\begin{remark}
In the case of hyperbolic surfaces, an important feature is that the bound is actually
independent of the volume, as long as the systole and spectral gap are controlled.

In the case of tori, it actually depends on the volume, as big volume flat tori must have
small Cheeger constant. 
\end{remark}
When dealing with zero average function, we cannot assume, of course, that they are nonnegative.
Hence we must be more careful when considering radial rearrangmeents:
\begin{theorem}\label{radialcompact}
	Let $(\Sigma,\sigma)$ be a closed surface.
	Debite $\delta,h,K$ its systole, Cheeger constant and supremum of curvature.
	Then there is a disk of radius $R=\frac{\Vol(\Sigma)}{2}$, with radial metric
	$dr^2+f'(r)^2d\theta^2$, so that for any $u$ Morse function on $\Sigma$ one can
	define $u_-$ and $u_+$ its lower and upper radial rearrangements satisfying:
	\begin{enumerate}
	\item
		$u_+$ is nonincreasing radial on $\D$.
	\item
		$u_-$ is nondecreasing radial on $\D$
	\item
	\begin{flalign}
		\forall t \text{ such that } \sigma(u\geq t)<\frac{\Vol(\Sigma)}{2},
		\sigma(u\geq t)&=|\{u_+\geq t\}|\\
		\forall t \text{ such that } \sigma(u\leq t)<\frac{\Vol(\Sigma)}{2},
		\sigma(u\leq t)&=|\{u_-\leq t\}|
	\end{flalign}
	\item
		For any measurable $f$,
	\begin{flalign}
		\int_\Sigma f(u)d\sigma&=\int_\D f(u_-)+f(u_+)\\
		\int_\Sigma |\nabla u|-2^2&\geq |\nabla u_+|_2^2+|\nabla u_-|_2^2
	\end{flalign}
	\end{enumerate}
\end{theorem}
\begin{proof}
	Denote $\varphi$ the isoperimetric profile of $\Sigma$.
	Thanks to lemma \ref{smallprofile}, there is $V_0>0$ so that,
	on $(0,V_0)$:
	\begin{equation}
		\varphi(t)\geq \sqrt{t(4\pi-K t)}
	\end{equation}
	Also, on $(0,\frac{\Vol(\Sigma)}{2}$, by definition of the Cheeger constant:
	\begin{equation}
		\varphi(t)\geq h\cdot t
	\end{equation}
	Consider the map $g$ defined in lemma \ref{globalprofile} with $\eps=V_0$.
	We check
	\begin{equation}
		\forall t\leq\frac{\Vol(\Sigma)}{2},\,\varphi(t)\geq g(t)
	\end{equation}
	Hence define $f$ by the lemma \ref{globalprofile}, and consider $u_+$ the nonincreasing
	radial rearrangement of the part of $u$ above its median $m$.
	In the same way, define $m-u_-$ to be the nonincreasing radial rearrangement of
	the part of $m-u$ above $u$.
	Thanks to the formulas \ref{formulas}, we deduce easily, for any measurable $f$:
	\begin{equation}
		\int_\Sigma f(u)=\int_\D f(u_+)+f(u_-)
	\end{equation}
	And because $g$ is lower than the isoperimetric profile on $(0,\frac{\Vol(\Sigma)}{2})$,
	we deduce:
	\begin{equation}
		\int_\Sigma |\nabla u|^2\geq |\nabla u_+|_2^2+|\nabla u_-|_2^2
	\end{equation}
\end{proof}
As a direct consequence, one can mimic the proof of the noncompact case and get:
\begin{corollary}
	Let $(\Sigma,\sigma)$ be a compact surface.
	Denote $\delta$, $h$ and $K$ its systole, Cheeger constant and supremum of curvature.
	There is a constant $C$ depending only on $\delta,h,K$ such that,
	for any $u\in W^{1,2}(\Sigma)$, satisfying:
	\begin{flalign}
		\int_\Sigma |\nabla u|^2&\leq 1\\
		\sigma(u\geq 0)&=\sigma(u\leq 0)
	\end{flalign}
	we have
	\begin{equation}
		\int_\Sigma \exp\big(4\pi u^2\big)-1\leq C
	\end{equation}
\end{corollary}
\begin{proof}
Consider $u$, and its radial rearrangements $u_+$,$u_-$ defined by theorem \ref{radialcompact}.
As the disk is conformal to a disk of curvature $-h$ with conformal bound depending only on
$r_0,h,K$, one can assume $\D$ is a disk with constant curvature $-h$.
We can then extend by zero the maps $u_+-m$ and $m-u_-$ to get functions in $W^{1,2}$ on the 
space of constant curvature $-h$.
Then , by \cite{MS10}, we have $C$ depending on $\delta, h, K$ such that:
\begin{flalign}
	\int_\D \exp\big(4\pi(u_+ -m)^2\big)-1&\leq C\\
	\int_\D \exp\big(4\pi(u_- -m)^2\big)-1&\leq C
\end{flalign}
The properties of the radial rearrangements are then enough to get the desired result:
\begin{equation}
	\int_\Sigma \exp\big(4\pi(u-m)^2\big)-1\leq C
\end{equation}
\end{proof}
In order to get the same statement for zero average functions,
we will need two technical lemmas.
\begin{lemma}\label{markov}
	Let $(\Sigma,\sigma)$ be a finite volume surface.
	Denote $\Lambda$ its spectral gap.
	Let $u\in W^{1,2}(\Sigma)$.
	Denote $\avg{u}$ its average and $m$ its median.
	Then
	\begin{equation}
		|m-\avg{u}|\leq\sqrt{\frac{2\Lambda}{\Vol(\Sigma)}}|\nabla u|_2
	\end{equation}
\end{lemma}
\begin{proof}
	Recall the Markov inequality, for $\eps>0$:
	\begin{equation}
	\sigma(\{|u-\avg{u}|\geq\eps\})\leq\eps^{-2}|u-\avg{u}|_2^2\leq\eps^{-2}\Lambda|\nabla u|_2^2
	\end{equation}	
	We apply it to
	\begin{equation}
		\eps^2=\frac{2\Lambda|\nabla u|_2^2}{\Vol(\Sigma)}
	\end{equation}
	in order to get
	\begin{equation}
		\sigma(\{|u-\avg{u}|\geq\eps\})\leq\frac{\Vol(\Sigma)}{2}\,.
	\end{equation}
	By definition of the median, necessarily
	\begin{equation}
		|\avg{u}-m|\leq\eps
	\end{equation}
	as asserted.
\end{proof}
\begin{lemma}\label{hardy}
	Let $\D$ be a disk of area $A$, and spectral gap $\Lambda$.
	Let $m>0$.
	Let $v$ be a radial nondecreasing map satisfying:
	\begin{flalign}
		v(\partial\D)&=m\\
		\int_\D v&\leq -m A\,. 
	\end{flalign}
	Then
	\begin{equation}
		|\nabla v|_2^2\geq 4\Lambda m^2 A
	\end{equation}
\end{lemma}
\begin{proof}
First, we use Cauchy-Schwarz for the estimate:
	\begin{equation}
		A^2 m^2\leq A|v|_2^2
	\end{equation}
	We deduce that
	\begin{equation}
		4A^2 m^2\leq|v-m|_2^2\,.
	\end{equation}
	Hence, by definition of the spectral gap:
	\begin{equation}
		|\nabla v|_2^2\geq 4\Lambda A^2 m^2
	\end{equation}
	as asserted.
\end{proof}
\begin{proof}[Proof of theorem \ref{TMcompact}]
	In order to get the result for zero average functions, we have to be more careful.
	Let $\Sigma$ be a compact surface,
	and consider $u$ a Morse function on it satisfying:
	\begin{flalign}
		\int_\Sigma u&=0\\
		\int_\Sigma |\nabla u|^2&\leq 1
	\end{flalign}
	Up to replacing $u$ by $-u$, we can assume its median satisfies:
	\begin{equation}
		m>0
	\end{equation}
	Consider its radial rearrangements $u_-$,$u_+$ on the disk $\D$,
	defined by theorem \ref{radialcompact}.
	As the disk is uniformly conformal to a hyperbolic disk,
	it has a spectral gap $\Lambda(\delta,h, K)$
	%one can assume that the disk has curvature $-h$.
	%Up to dilating by $h$, whose effect will be to multiply the upper bound by
	%$\frac{1}{h^2}$, we can assume $\D$ to be a hyperbolic disk of radius $A$.
	As $u$ is of zero average, $u_-$ satisfies the assumptions of lemma \ref{hardy},
	and we have:
	\begin{equation}
		|\nabla u_-|_2^2\geq 4\Lambda m^2 A
	\end{equation}
	where $A$ denotes the area of our disk.
	Of course, this implies
	\begin{equation}
		|\nabla u_+|_2^2\leq 1- 4\Lambda m^2 A
	\end{equation}
	We can now apply the Moser-Trudinger inequality to
	$\frac{1}{\sqrt{1-m^2A}}(u_+-m)$ to get
	\begin{equation}
		\int_\D\exp\big(\frac{4\pi}{1-4\Lambda m^2 A}(u_+-m)^2\big)-1\leq C(\delta, h, K)
	\end{equation}
	We compute
	\begin{flalign}
		\frac{(u-m)^2}{1-4 \Lambda m^2 A}-u^2&=\frac{4\Lambda m^2 A}{1-4\Lambda m^2 A}\bigg[u^2-\frac{2u}{4\Lambda mA}+\frac{1}{4\Lambda A}\bigg]\\
		&\geq-\frac{1}{4\Lambda A}
	\end{flalign}
	and we deduce that
	\begin{equation}
		\int_\D\exp\big(4\pi u_+^2\big)-1\leq e^{\frac{1}{4\Lambda A}}C
		+4\Lambda A(e^{\frac{1}{4\Lambda A}}-1)
	\end{equation}
	Byy construction,
	\begin{equation}
		2A=\Vol(\Sigma)
	\end{equation}
	Hence we get a bound $C'$ depending on $(\delta, h,K)$ only so that:
	\begin{equation}
		\int_\D\exp\big(4\pi u_+^2\big)-1\leq C'\exp\big(\frac{1}{2\Lambda\Vol(\Sigma)}\big)
			+2\Lambda\Vol(\Sigma)\big(\exp\frac{1}{2\Lambda\Vol(\Sigma)}-1\big)
	\end{equation}
	Finally, we check that
	\begin{equation}
		u_-^2\leq (u_- -m)^2+m^2
	\end{equation}
	And by lemma \ref{markov},
	\begin{equation}
		m\leq \sqrt{\frac{2\Lambda}{\Vol(\Sigma)}}\,.
	\end{equation}
	Hence we can bound
	\begin{flalign}
		\int_\D\exp\big(4\pi u_-^2\big)-1&\leq C(\delta,h,K)e^{m^2}
		+\frac{\Vol(\Sigma)}{2}(e^{m^2}-1)\\
		&\leq C(\delta,h,K)e^\frac{2\Lambda}{\Vol(\Sigma)}
		+\frac{\Vol(\Sigma)}{2}(e^\frac{2\Lambda}{\Vol(\Sigma)}-1)
	\end{flalign}
	All in all, we get nonnegative constants $C'',a,b$ depending on $\delta, h, K$ so that:
	\begin{equation}
		\int_\D\exp\big(4\pi u^2\big)-1\leq 
		C''(\exp\frac{a}{\Vol(\Sigma)})+\Vol(\Sigma)(\exp\frac{b}{\Vol(\Sigma)}-1)
	\end{equation}
	The map $t\mapsto t(\exp\frac{b}{t}-1)$ is nonincreasing on $\R_+-\{0\}$,
	hence we get an upper bound which is nonincreasing in $\Vol(\Sigma)$, as asserted.
\end{proof}

We can see that these dependencies are actually sharp for the case of hyperbolic closed surfaces:
\begin{proposition}
	Let $\eps_n\rightarrow 0$ and let $(X_n)$ be a family of hyperbolic closed surfaces
	with systole $\eps_n$.
	Then there is a family of functions $f_n\in W^{1,2}(X_n)$ satisfying:
	\begin{equation}
		\int_{X_n} f_n=0,\quad\int_{X_n}|\nabla f_n|^2\leq 1
	\end{equation}
	And
	\begin{equation}
		\int_{X_n}|f_n|^4\rightarrow +\infty
	\end{equation}
\end{proposition}
\begin{proof}
	Fix $R_0>0$.
	Thanks to the collar lemma (see \cite{Bus10}, theorem 4.1.1), if the systole goes to zero,
	for $\eps$ small enough $X_n$ contains an isometric collar
	$(-R_0,R_0)\times\S^1$ with metric $dt^2+\eps_n^2ch(t)^2d\theta^2$.
	Consider the map:
	\begin{equation}
		g_n(t,\theta)=\one_{(-R_0,R_0)}(\cosh(t)-\cosh(R_0)
	\end{equation}
	It belongs to $W^{1,2}(X_n)$, and satisfies:
	\begin{flalign}
		\int_{X_n}|\nabla g_n|^2&=\frac{2}{3}\eps_n\cosh(R_0)\\
		\int_{X_n} g_n&=2\eps_n(\frac{-\sinh(2R_0)}{4}+\frac{R_0}{2})\\
		\int_{X_n} |g_n|^2&=2\eps_n(\frac{\sinh(R_0)^3}{3}+\sinh(R_0)-R_0\cosh(R_0))
	\end{flalign}
	A thorough computation gives, when $\eps_n\rightarrow 0$ and $R_0$ too, with $\eps_n$ negligible with regard to $R_0$:
	\begin{equation}
		\int_{X_n} |g_n-\frac{1}{\Vol(\Sigma)}\int_{X_n} f|^4\sim \frac{16\eps_n R_0^9}{315}
	\end{equation}
	All in all, fix $R_0=\eps^\frac{1}{10}$ small enough and you get, as $n\rightarrow +\infty$:
	\begin{equation}
		\frac{|g_n-\avg{g_n}|_4^4}{|\nabla g_n|_2^4}\sim \frac{4}{35}\eps^{-\frac{1}{10}}\rightarrow+\infty
	\end{equation}
	as asserted.
\end{proof}
In particular, if we consider a family with lower bounded volume and systole,
and uniformly bounded curvature and Cheeger constant, we get a Moser-Trudinger bound independent
of the surface considered in the family.

\begin{corollary}
	Let $X$ be a closed hyperbolic surface, and $\eps>0$
	Consider $\mathcal{X}$ the collection of Riemannian covers of $X$ with Cheeger constant
	bigger than $\eps$.
	Then there is a constant $C(\eps)$ such that,
	for any $Y\in\mathcal{X}$, for any $u\in W^{1,2}(Y)$ satisfying:
	\begin{equation}
		\int_Y u=0\quad\int_Y|\nabla u|^2\leq 1
	\end{equation}
	we have
	\begin{equation}
		\int_Y(e^{4\pi u^2}-1)\leq C(\eps)
	\end{equation}
\end{corollary}

\begin{remark}
Now depending on $\eps$ this family has different characteristics:
For $\eps>h_1(X)$, this family is empty.
The result of \cite{MNP22} combined with Buser's control tells us that for $\eps$
close enough to $0$,
the family is infinite, and even "big", in the sense that the probability that a random cover of 
degree $d$ doesn't belong to $\mathcal{X}$ goes to zero as $d$ is large.

It is not known to the author if a random tower of covers of $X$ has a nonzero probability to
belong to some $\mathcal{X}(\eps)$, for $\eps>0$.
\end{remark}

\newpage
\bibliographystyle{alpha}
\bibliography{references}

\newcommand{\etalchar}[1]{$^{#1}$}
\begin{thebibliography}{BJMR15}

\bibitem[BJMR15]{BJM+15}
Luca Battaglia, Aleks Jevnikar, Andrea Malchiodi, and David Ruiz.
\newblock A general existence result for the toda system on compact surfaces.
\newblock {\em Advances in Mathematics}, 285:937--979, 2015.

\bibitem[Bor07]{Bor07}
David Borthwick.
\newblock {\em Spectral theory of infinite-area hyperbolic surfaces}.
\newblock Springer, 2007.

\bibitem[BPP{\etalchar{+}}93]{BPP93}
Robert Brooks, Peter Perry, Peter Petersen, et~al.
\newblock On cheeger’s inequality.
\newblock {\em Comment. Math. Helv}, 68(4):599--621, 1993.

\bibitem[Bus79]{Bus79}
Peter Buser.
\newblock Ueber den ersten eigenwert des laplace-operators auf kompakten
  fl{\"a}chen.
\newblock {\em Commentarii Mathematici Helvetici}, 54(1):477--493, 1979.

\bibitem[Bus82]{Bus82}
Peter Buser.
\newblock A note on the isoperimetric constant.
\newblock In {\em Annales scientifiques de l'{\'E}cole normale sup{\'e}rieure},
  volume~15, pages 213--230, 1982.

\bibitem[Bus10]{Bus10}
Peter Buser.
\newblock {\em Geometry and spectra of compact Riemann surfaces}.
\newblock Springer Science \& Business Media, 2010.

\bibitem[Che70]{Che70}
J~Cheeger.
\newblock A lower bound for the smallest eigenvalue of the laplacian.
\newblock {\em Problems in Analysis}, pages 195--199, 1970.

\bibitem[DEJ14]{DEJ14}
Jean Dolbeault, Maria~J Esteban, and Gaspard Jankowiak.
\newblock Onofri inequalities and rigidity results.
\newblock {\em arXiv preprint arXiv:1404.7338}, 2014.

\bibitem[Fag06]{Fag06}
Zo{\'e} Faget.
\newblock Best constants in the exceptional case of sobolev inequalities.
\newblock {\em Mathematische Zeitschrift}, 252(1):133--146, 2006.

\bibitem[Gal88]{Gal88}
Sylvestre Gallot.
\newblock In{\'e}galit{\'e}s isop{\'e}rim{\'e}triques et analytiques sur les
  vari{\'e}t{\'e}s riemanniennes.
\newblock {\em Ast{\'e}risque}, 163(164):31--91, 1988.

\bibitem[GNP09]{GNP09}
Renata Grimaldi, Stefano Nardulli, and Pierre Pansu.
\newblock Semianalyticity of isoperimetric profiles.
\newblock {\em Differential Geometry and its Applications}, 27(3):393--398,
  2009.

\bibitem[Ili83]{Ili83}
Sa{\"\i}d Ilias.
\newblock Constantes explicites pour les in{\'e}galit{\'e}s de sobolev sur les
  vari{\'e}t{\'e}s riemanniennes compactes.
\newblock In {\em Annales de l'institut Fourier}, volume~33, pages 151--165,
  1983.

\bibitem[JW01]{JW01}
J{\"u}rgen Jost and Guofang Wang.
\newblock Analytic aspects of the toda system: I. a moser-trudinger inequality.
\newblock {\em Communications on Pure and Applied Mathematics: A Journal Issued
  by the Courant Institute of Mathematical Sciences}, 54(11):1289--1319, 2001.

\bibitem[Kes06]{Kes06}
Srinivasan Kesavan.
\newblock {\em Symmetrization and applications}, volume~3.
\newblock World Scientific, 2006.

\bibitem[Li01]{Li01}
Yuxiang Li.
\newblock Moser-trudinger inequality on compact riemannian manifolds of
  dimension two.
\newblock {\em Journal of Partial Differential Equations}, 14(2):163--192,
  2001.

\bibitem[MNP22]{MNP22}
Michael Magee, Fr{\'e}d{\'e}ric Naud, and Doron Puder.
\newblock A random cover of a compact hyperbolic surface has relative spectral
  gap 3 16-$\varepsilon$.
\newblock {\em Geometric and Functional Analysis}, 32(3):595--661, 2022.

\bibitem[Mos71]{Mos71}
J{\"u}rgen Moser.
\newblock A sharp form of an inequality by n. trudinger.
\newblock {\em Indiana University Mathematics Journal}, 20(11):1077--1092,
  1971.

\bibitem[Mos73]{Mos73}
J{\"u}rgen Moser.
\newblock On a nonlinear problem in differential geometry.
\newblock In {\em Dynamical systems}, pages 273--280. Elsevier, 1973.

\bibitem[MR11]{MR11}
Andrea Malchiodi and David Ruiz.
\newblock A variational analysis of the toda system on compact surfaces.
\newblock {\em arXiv preprint arXiv:1105.3701}, 2011.

\bibitem[MS10]{MS10}
G~Mancini and K~Sandeep.
\newblock Moser--trudinger inequality on conformal discs.
\newblock {\em Communications in Contemporary Mathematics}, 12(06):1055--1068,
  2010.

\bibitem[Tru67]{Tru67}
Neil~S. Trudinger.
\newblock On imbeddings into orlicz spaces and some applications.
\newblock {\em Journal of Mathematics and Mechanics}, 17(5):473--483, 1967.

\bibitem[Xio18]{Xio18}
Jingang Xiong.
\newblock A derivation of the sharp moser--trudinger--onofri inequalities from
  the fractional sobolev inequalities.
\newblock {\em Peking Mathematical Journal}, 1(2):221--229, 2018.

\end{thebibliography}
\end{document}